\documentclass{amsart}
\usepackage{amssymb}
\usepackage[cp1250]{inputenc}
\usepackage{graphicx}
\newtheorem{thm}{Theorem}[section]

\newtheorem{lem}[thm]{Lemma}

\theoremstyle{definition}

\theoremstyle{remark}
\newtheorem{rem}[thm]{Remark}
\numberwithin{equation}{section}

\newcommand{\Real}{\mathbb R}
\newcommand{\Integer}{\mathbb Z}

\newcommand{\To}{\longrightarrow}

\newcommand{\DM}{\begin{displaymath}}
\newcommand{\EDM}{\end{displaymath}}
\newcommand{\BE}{\begin{equation}}
\newcommand{\EE}{\end{equation}}

\newcommand{\comment}[1]{\mbox{}}

\newcommand{\diam }{\,{\rm diam}\,}

\begin{document}
\title[Simple choreographies in N-body problem]
 {An existence of simple choreographies for N-body problem - a computer assisted proof}

\author{Tomasz Kapela}
\author{Piotr Zgliczy\'nski}

\address{Tomasz Kapela, Pedagogical University, Institute of Mathematics,
\mbox{Podchora\.zych 2}, 30-084 Krak\'ow, Poland }
\email{tkapela@ap.krakow.pl}

\address{Piotr Zgliczy\'nski, Jagiellonian University, Institute of Mathematics, Reymonta
4, 30-059 Krak\'ow, Poland}
\email{zgliczyn@im.uj.edu.pl}

\thanks{}



\date{\today}

\begin{abstract}
We consider a question of finding a periodic solution for the
planar Newtonian N-body problem with equal masses, where each body
is travelling along the same closed path. We provide a computer
assisted proof for the following facts: local uniqueness and
convexity of Chenciner and Montgomery Eight, an existence (and
local uniqueness) for Gerver's SuperEight for 4-bodies and a
doubly symmetric linear chain for 6-bodies.
\end{abstract}

\maketitle

\section{Introduction}
\label{sec:intro}  In this paper we consider the problem of
finding  periodic solution to the $N$-body problem in which all
$N$ masses travel along a fixed curve in the plane. The $N$-body
problem with $N$ equal unit masses is given by a differential
equation
\begin{equation}
 \ddot{q_i}=\sum_{j\neq i} \frac{q_j -q_i}{r_{ij}^3}\label{eq:nbp}
\end{equation}
where $q_i\in \Real^n$, $i=1,...,N$, $r_{ij}= ||q_i - q_j||$. The
gravitational constant is taken equal to 1.

\par
We  consider planar case (n=2) only,  we set $q_i=(x_i, y_i)$,
$\dot{q_i}=(v_i, u_i)$. Using this we can express (\ref{eq:nbp})
by: \BE \left\{\begin{array}{l}
\dot{v_i}=\sum_{j\neq i}\frac{x_j -x_i}{r_{ij}^3}\\
\dot{u_i}=\sum_{j\neq i}\frac{y_j -y_i}{r_{ij}^3}\\
\dot{x_i}=v_i\\
\dot{y_i}=u_i
\end{array} \right.
\EE Recently, this problem received a lot attention in literature
see \cite{M,CM,CGMS,S1,S2,S3,MR} and papers cited there.

 By a {\em simple choreography} \cite{S1,S2} we mean a collision-free
 solution of the N-body
problem in which all masses move on the same curve with a constant
phase shift. This means that there exists $q:\Real \To \Real^2$ a
$T$-periodic function of time, such that the position of $k$-th
body ($k = 0,\dots,N-1$) is given by $q_k(t) = q(t + k
\frac{T}{N})$ and $(q_0, q_1, \dots, q_{N-1})$ is solution of the
$N$-body problem. The simplest choreographies are  Langrange
solutions in which the bodies are located at the vertices of a
regular $N$-gon and move with constant angular velocity. Another
simple choreography, a figure eight curve (see
Figure~\ref{fig:eight}), was found numerically by C.Moore
\cite{M}.  A.Chenciner and R.Montgomery \cite{CM} gave a rigorous
existence proof of the Eight in 2000. In December 1999, J.Gerver
found orbit for $N=4$ called 'Super-Eight'
(Figure~\ref{fig:gerver}). After that C.Sim\'o found a lot of
simple choreographies with different shape and number of bodies
ranging form 4 to several hundreds (see \cite{S,S1,CGMS} for
pictures, animations and more details).

Up to now the only choreographies whose existence has been
rigorously established \cite{MR} are Lagrange solutions and the
Eight solution. While the Lagrange solution is given analytically,
the existence of Eight was proven in \cite{CM} using variational
arguments and  there is still a lot of open questions about it
\cite{Ch,MR}. For example uniqueness (up to obvious symmetries and
rescaling ) and convexity of the lobes in Eight. In
Section~\ref{sec:eight} we give a computer assisted proof of
existence of Eight, its local uniqueness and  convexity of the
lobes.

In Sections \ref{sec:linchains}, \ref{sec:gerver} and
\ref{sec:6body} we concentrate on choreographies called {\em
doubly symmetric linear chains}. They are symmetric with respect
to the $x$ and $y$ coordinate axes and all points of the curve
self-intersection are on the $x$ axis.
 We give a computer assisted proof of an existence (and local uniqueness)
 of doubly symmetric linear chains for four  (Gerver SuperEight) and six bodies.

Proofs given in this paper are computer assisted. By this we mean
that we use computer to provide rigorous bounds for solutions of
(\ref{eq:nbp}). The problem of proving an existence of a
choreography is reduced to finding a zero for a suitable function.
For this purpose we use an interval Newton method and a Krawczyk
method \cite{A,K,Mo,N}, (see Section~\ref{sec:inter-thm}) which
are apparently rather unknown outside the interval arithmetic
community. To integrate equations (\ref{eq:nbp}) we use a
$C^1$-Lohner algorithm \cite{ZLo}. All computations were performed
on {\em AMD Athlon 1700XP } with 256 MB DDRAM memory, with Windows
98SE operating system. We used CAPD package \cite{Capd} and
Borland C++ 5.02 compiler. A total computation time for 6-bodies
was under 90 seconds and considerably smaller for the Eight and
SuperEight solutions - see Section~\ref{sec:tech-data} for more
details.

\section{Two zero finding methods}
\label{sec:inter-thm} The main technical tool used in this paper
in order to establish an existence of solutions of equations of
the form $f(x)=0$ is {\em an interval Newton method}\cite{A,Mo,N}
and {\em a Krawczyk method}\cite{A,K,N}. The interval Newton
method was used to prove an existence of the Eight
(Figure~\ref{fig:eight}) and Gerver orbit
(Figure~\ref{fig:gerver}). The Krawczyk method for Gerver orbit
and an orbit with 6 bodies in a linear chain
(Figure~\ref{fig:6body}).

\subsection{Notation}
In the application of  interval arithmetics to rigorous
verification of theorems single valued objects, like numbers,
vectors, matrices etc are in the formulas replaced by sets
containing sure bounds for them.
 In the sequel, we will not use any special notation for single
valued object and sets. For a set $S$ by $[S]$ we denote the
interval hull of $S$, i.e. the smallest product of intervals
containing $S$. For a set which is an interval set (i.e. can be
represented as a product of intervals) we will also use square
brackets to stress its interval nature. For any interval set $[S]$
by $mid([S])$ we denote a center point of $[S]$. For any interval
$[a,b]$ we define a diameter by $diam [a,b] := b - a $. For an
interval vector (matrix) $S=[S]$ by $diam S$ we denote a vector
(matrix) of diameters of each components.

\subsection{Interval Newton method}

\begin{thm} \cite{A,N}
\label{thm:newton} Let $F:\Real^n \to \Real^n$ be a $C^1$
function. Let $[X]=\Pi_{i=1}^n[a_i,b_i]$, $a_i < b_i$. Assume the
interval hull of $DF([X])$, denoted here by $[DF([X])]$, is
invertible. Let $\bar{x} \in X$ and we define
\begin{equation}
  N(\bar{x},[X]) = - [DF([X])]^{-1}F(\bar{x}) + \bar{x}
\end{equation}
Then
\begin{enumerate}
  \item[0.] if \ $x_1,x_2 \in [X]$ and $F(x_1)=F(x_2)$, then $x_1=x_2$
  \item[1.] if \ $N(\bar{x},[X]) \subset [X]$, then $\exists ! x^* \in [X]$ such that
   $F(x^*)=0$
  \item[2.] if \ $x_1 \in [X]$ and $F(x_1)=0$, then $x_1 \in N(\bar{x},[X])$
  \item[3.] if \ $N(\bar{x},[X]) \cap [X] = \emptyset$, then $F(x) \neq 0$ for all $x \in [X]$
\end{enumerate}
\end{thm}

\subsection{Krawczyk method}
We assume that:
\begin{itemize}
  \item $F:\Real^n \To \Real^n$ is a $C^1$ function,
  \item $[X] \subset \Real^n$ is an interval set,
  \item $\bar{x} \in [X]$
  \item $C \in \Real^{n\times n}$ is a linear isomorphism.
\end{itemize}
The Krawczyk \cite{A,K,N} operator is given by
\begin{equation}\label{eq:krawczyk_operator}
  K(\bar{x}, [X], F) := \bar{x} - CF(\bar{x}) + (Id -
  C[DF([X])])([X]-\bar{x}).
\end{equation}
\begin{thm}\label{thm:Krawczyk}
1. If $x^* \in [X]$ and $F(x^*) = 0$, then $x^* \in K(\bar{x},
[X], F)$.

2. If $K(\bar{x}, [X], F) \subset int [X]$, then there exists in
[X] exactly one solution of equation $F(x)=0$.

3. If $K(\bar{x}, [X], F) \cap  [X]=\emptyset$, then $F(x) \neq 0$
for all $x \in [X]$
\end{thm}

\subsection{Algorithm for Newton and Krawczyk method}
Theorems~\ref{thm:newton} and \ref{thm:Krawczyk} can be used as a
basis for an algorithm for rigorous enclosing for solution of
equation $F(x)=0$. Let $T(\bar{x},[X])=N(\bar{x},[X])$ if we are
using the interval Newton method and $
T(\bar{x},[X])=K(\bar{x},[X],F)$ in case of Krawczyk method.

 First, we need to have a good guess for ${x}^* \in \Real^n$.
For this purpose we   use a nonrigorous Newton method to obtain
$\bar{x}=(\bar{x}_1,\bar{x}_2,\dots,\bar{x}_n)$. Then we choose
interval set $[X]$ which contains $\bar{x}$ and perform the
following algorithm:
\begin{itemize}
\item[Step 1.] Compute $T(\bar{x},[X])$.
\item[Step 2.] If $T(\bar{x},[X]) \subset [X]$, then return
\textbf{success}.
\item[Step 3.] If $X \cap T(\bar{x},[X]) = \emptyset$, then return
\textbf{fail}. There are no zeroes of $F$ in $[X]$.
\item[Step 4.] If $[X] \subset T(\bar{x},[X]) $, then modify
computation parameters (for example: a time step, the order of
Taylor method, size of $[X]$). Go to Step 1.
\item[Step 5.] Define a new $[X]$ by $[X] := [X] \cap
T(\bar{x},[X])$ and a new $\bar{x}$ by $\bar{x} := mid ([X])$,
then go to Step 1.
\end{itemize}

In practical computation it is convenient to define a maximum
number of iteration allowed and return \textbf{fail} if the actual
iteration count is larger.

Observe that the third assertion in both Theorems \ref{thm:newton}
and \ref{thm:Krawczyk} can be used to exclude an existence of zero
of $F$. This has been used by Galias  \cite{Gp1,Gp2} to find all
periodic orbits up to a given period for H\'enon map and Ikeda
map. It seems possible to obtain similar results for N-body
problem in the future.

\section{The Eight - existence, local uniqueness and convexity}
\label{sec:eight} Existence of the Eight has been shown in
\cite{CM} using mixture of symmetry and variational arguments.
Here we give an another existence proof and in addition we obtain
local uniqueness and the convexity of each lobe of the Eight. We
follow \cite{CM} in the use of the symmetry, but other component
of the proof is different - we use an interval Newton method
discussed in Section~\ref{sec:inter-thm}.
\par
In notation, but only for the Eight, we follow \cite{CM}.

\begin{figure} [hbdt]
 \centerline{\includegraphics[scale = 0.7]{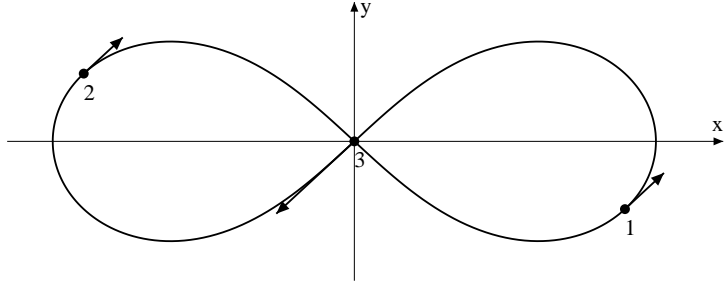}}
 \caption{The Eight - the initial position \label{fig:eight}}
\end{figure}
\begin{figure}  

 \centerline{ \includegraphics[scale = 0.7]{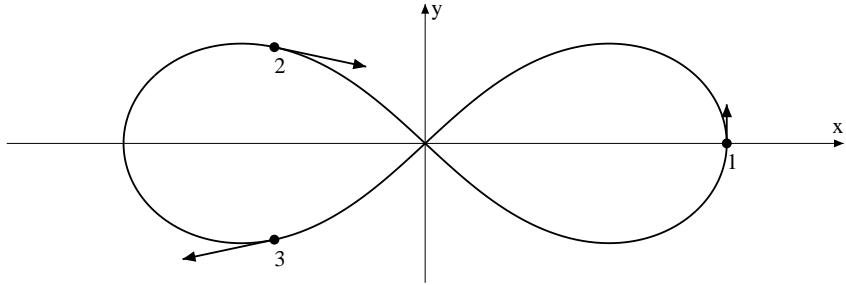}}
 \caption{The Eight - the final position \label{fig:eight_e}}
\end{figure}

Let $T$ be any positive real number. We define action of the Klein
group $\Integer_2 \times \Integer_2$ on $\Real/T\Integer$ and on
$\Real^2$ as follows: if $\sigma$ and $\tau$ are generators,
\begin{equation}
\sigma(t) = t + \frac{T}{2},\  \tau(t) = -t + \frac{T}{2},\
\sigma(x, y) = (-x, y),\  \tau(x,y) = (x, -y). \label{eq:8-sym}
\end{equation}

For loop $q:(\Real/T\Integer) \To \Real^2$ and $i = 1,2,3$ we
define the position for $i$-th body by
\begin{equation}
  q_i(t) = q\left(t + (3-i) \cdot \frac{T}{3}\right),
\end{equation}
Here $q_i$ is just position of $i$-th body.

The following theorem without the uniqueness part was proved in
\cite{CM}.
\begin{thm}
\label{thm:eight} There exists an "eight"-shaped planar loop
$q:(\Real/T\Integer) \To \Real^2$ with the following properties:
\begin{enumerate}
\item
for each t,
\[
q_1(t) + q_2(t) + q_3(t) = 0;
\]
\item
q is invariant with respect to the action of $\Integer_2 \times
\Integer_2$ on $\Real/T\Integer$ and on $\Real^2$:
\[
q \circ \sigma(t) = \sigma \circ q(t) \mbox{ and } q \circ
\tau(t)=\tau \circ q(t);
\]
\item
the loop $x: \Real/T\Integer \To \Real^6 $ defined by
\[
x(t) = (q_1(t), q_2(t), q_3(t))
\]
is $T$-periodic solution of the planar three-body problem with
equal masses.

Moreover, $q$ is locally unique (up to obvious symmetries  and
rescaling).
\end{enumerate}
\end{thm}

\begin{rem}
\label{rem:eight} If conditions (1),(2),(3) are satisfied then
\begin{enumerate}
\item[a.] $\dot q_1(t) + \dot q_2(t) + \dot q_3(t) = 0$ for each t,
\item[b.] $\dot{q} \circ \sigma(t) = \sigma \circ \dot{q}(t) \mbox{ and }
\dot q  \circ \tau(t)=\sigma \circ \dot{q}(t)$ for each t,
\item[c.] $q_3(0)=(0,0)$ and $\dot{q_3}(0)=-2\dot{q_1}(0)$
\item[d.] $q_1(0) = - q_2(0)$ and $\dot q_1(0) = \dot q_2(0)$,
\item[e.] $q_1(T/12)$ is on the X axis and $\dot q_1(T/12)$ is orthogonal to the X axis,
\item[f.] $q_3(T/12) = \tau \circ q_2(T/12)$ and $\dot q_3(T/12) = \sigma \circ \dot q_2(T/12)$,
\end{enumerate}
\end{rem}

The following lemma describes the symmetry reduction for the Eight
\begin{lem}
\label{lem:eight} Assume that $\tilde{q}: [0, \tilde{T}] \To
\Real^6$ is a solution of the three body problem, such that
$\tilde{q}=(\tilde{q}_1, \tilde{q}_2, \tilde{q}_3)$ satisfies
conditions (c),(d),(e),(f) in Remark~\ref{rem:eight}, then exists
$q:(\Real/T\Integer)\To \Real^2$ satisfying condition (1),(2),(3)
in Theorem~\ref{thm:eight} with $T=12\tilde{T}$.
\end{lem}
\noindent \textbf{Proof:} We define
\begin{equation}\label{def:q_1_4}
\hat{q}(t)=\left\{ \begin{array}{ll}
  \tilde{q}_3(t) & \mbox{for } t \in [0,\tilde{T}] \\
  \tau \circ \tilde{q}_2(2\tilde{T} - t) & \mbox{for } t \in [\tilde{T},2\tilde{T}] \\
  \sigma \circ \tilde{q}_1(t-2\tilde{T}) & \mbox{for } t \in [2\tilde{T},3\tilde{T}] \\
\end{array} \right.
\end{equation}
and
\begin{equation}\label{def:q}
q(t)=\left\{ \begin{array}{ll}
  \hat{q}(t) & \mbox{for } t \in [0,3\tilde{T}] \\
  \tau \circ \hat{q}(\tau^{-1}(t)) & \mbox{for } t \in [3\tilde{T},6\tilde{T}] \\
  \sigma \circ \hat{q}(\sigma^{-1}(t)) & \mbox{for } t \in [6\tilde{T},9\tilde{T}] \\
  \sigma \circ \tau \circ \hat{q}(\tau^{-1} \circ \sigma^{-1}(t))& \mbox{for } t \in [9\tilde{T},12\tilde{T}] \\
\end{array} \right.
\end{equation}

Let $f_1$ and $f_2$ be two solutions of \ref{eq:nbp} on intervals
$[t_1, t_2]$ and $[t_2, t_3]$ respectively. If $f_1(t_2) =
f_2(t_2)$ and $\dot{f}_1^-(t_2) = \dot{f}_2^+(t_2)$ then
$f=\{f_1,f_2\}$ is solution on interval $[t_1,t_2]$. To show that
$q(t)$ is a solution it is enough to show that "pieces fit
together smoothly".

\noindent For $t=\tilde{T}$ from \ref{rem:eight}.f we have
\begin{eqnarray}
\label{eqn:qt1}
 \tilde{q}_3(\tilde{T}) &=& \tilde{q}_3(\frac{T}{12})= \tau \circ
 \tilde{q}_2(\frac{T}{12})=\tau \circ \tilde{q}_2(\tilde{T})=\tau
 \circ \tilde{q}_2(2\tilde{T}-\tilde{T}) \\
 \dot{\tilde{q}}_3^-(\tilde{T}) &=& \sigma \circ
 \dot{\tilde{q}}_2^+(\tilde{T})
\end{eqnarray}

\noindent For $t=2\tilde{T}$ from \ref{rem:eight}.d we have
\begin{eqnarray}
 \tilde{q}(2\tilde{T})&=& \tau \circ \tilde{q}_2(0) =
 \sigma \circ \tilde{q}_1(0)  \\
 \label{eqn:qt4}
 \dot{\tilde{q}}^-(2\tilde{T}) &=& \sigma \circ \dot{\tilde{q}}_2^-(\tilde{T}) = \sigma \circ
 \dot{\tilde{q}}_1^+(2\tilde{T}) =  \dot{\tilde{q}}^+(2\tilde{T})
\end{eqnarray}

From (\ref{eqn:qt1}-\ref{eqn:qt4}) it follows  that $\hat{q}(t)$
is smooth curve for $t \in [0, 3\tilde{T}]$. We construct $q(t)$
from that curve.

\noindent For $t=3\tilde{T}$ from \ref{rem:eight}.e we obtain
\begin{eqnarray*}
 {q}(3\tilde{T})&=& \sigma \circ \tilde{q}_1(\tilde{T}) = \sigma \circ \tau \circ \tilde{q}_1(\tilde{T}) =
 \sigma \circ \tau \circ \hat{q}(\tau^{-1}(3\tilde{T}))  \\
 \dot{q}^-(3\tilde{T}) &=& \sigma \circ \dot{\tilde{q}}_1(\tilde{T}) =
 \dot{\tilde{q}}_1(\tilde{T}) =   \dot{q}^+(3\tilde{T})
\end{eqnarray*}

\noindent Using \ref{rem:eight}.c for $t=6\tilde{T}$  we infer
that
\begin{eqnarray*}
 {q}(6\tilde{T})&=& \tau \circ \hat{q}(\tau^{-1}(6\tilde{T})) = \tau \circ \tilde{q}_3(0) =
 \sigma \circ \tilde{q}_3(0) =  \sigma \circ
 \hat{q}(\sigma^{-1}(6\tilde{T})) \\
 \dot{q}^-(6\tilde{T}) &=& \sigma \circ \dot{\tilde{q}}_1(\tilde{T}) =  \dot{q}^+(6\tilde{T})
\end{eqnarray*}

\noindent Now we show that $q(t)$ is a closed curve.
\begin{eqnarray*}
 {q}(12\tilde{T})= \sigma \circ \tau \circ \hat{q}(\tau^{-1}\circ\sigma^{-1}(12\tilde{T}))
 = \sigma \circ \tau \circ \hat{q}(0) = \hat{q}(0) = q(0)\\
 \dot{q}^-(12\tilde{T}) = \frac{\partial}{\partial t}[\sigma \circ \tau \circ
 \hat{q}(\tau^{-1}\circ\sigma^{-1}(12\tilde{T}))]=
 \dot{\hat{q}}(0) = \dot{q}^+(0)
\end{eqnarray*}
Hence it is easy to see that $q(t)$ can be extended to a
$T$-periodic curve, such that $x(t) = (q_1(t), q_2(t), q_3(t))$ is
a $T$-periodic solution of the three body problem.

Condition (2) in Theorem~\ref{thm:eight} follows easily from
definition of $\hat{q}(t)$ and properties of $\sigma$ and $\tau$
(for T-periodic orbit $\sigma^{-1}(t)=\sigma(t)$ and
$\tau^{-1}(t)=\tau(t)$). For $t \in [0, 3\tilde{T}]$ we have
\begin{eqnarray}
q(\sigma(t)) &=& \sigma \circ \hat{q}(\sigma^{-1} \circ \sigma(t))
=
\sigma \circ \hat{q}(t) = \sigma \circ q(t) \\
\label{eqn:w2_2} q(\tau (t)) &=& \tau \circ \hat{q}(\tau^{-1}\circ
\tau(t)) = \tau \circ \hat{q}(t) = \tau \circ q(t)
\end{eqnarray}
For $t \in [3\tilde{T}, 6\tilde{T}]$ using (\ref{eqn:w2_2}) we
obtain
\begin{eqnarray}
q(\sigma(t)) = \sigma \circ \tau \circ \hat{q}(\tau^{-1} \circ
\sigma^{-1} \circ \sigma(t)) = \sigma \circ \tau \circ
\hat{q}(\tau^{-1}(t)) = \sigma \circ q(t) \\
q(\tau (t)) = \hat{q}(\tau(t)) = \tau \circ \hat{q}(\tau \circ
\tau (t)) = \tau \circ q(t)
\end{eqnarray}
We omit other two cases, because the proof is very similar.

Because there aren't any external force so  the center of mass can
only move with a constant velocity. But from \ref{rem:eight}.c and
\ref{rem:eight}.d we have $q_1(0) + q_2(0) + q_3(0) = 0$ and $\dot
q_1(0) + \dot q_2(0) + \dot q_3(0) = 0$. Hence we obtain condition
(1) in Theorem~\ref{thm:eight}.
 \qed

\begin{figure}[hdtb]
 \centerline{\includegraphics[scale = 0.7]{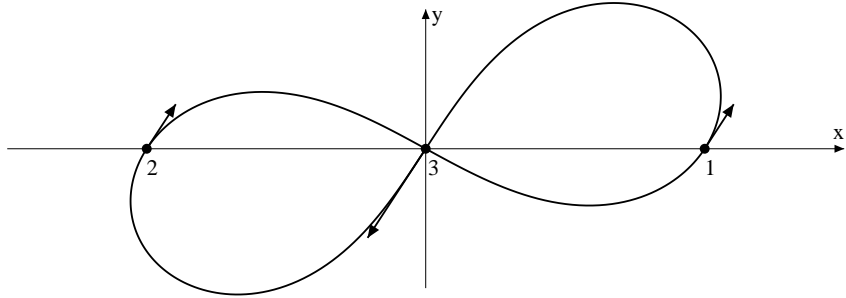}}
 \caption{Rotated Eight - initial position \label{fig:eight2}}
\end{figure}
\begin{figure}[hdtb]
 \centerline{\includegraphics[scale = 0.7]{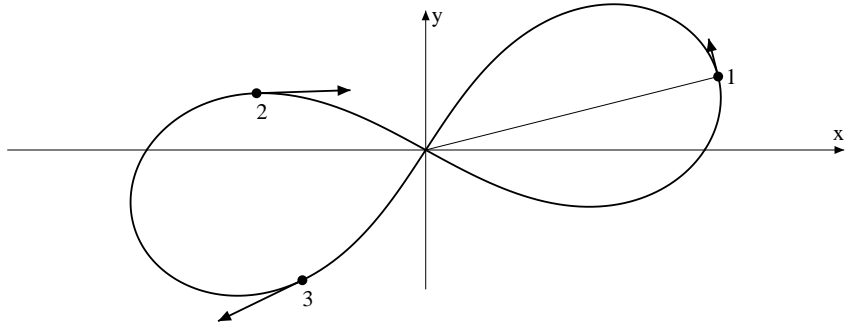}}
 \caption{Rotated Eight - final position \label{fig:eight2_e}}
\end{figure}

Hence to prove Theorem~\ref{thm:eight} it is enough to show that
there exists a locally unique (up to obvious degeneracies)
function satisfying assumptions of Lemma~\ref{lem:eight}. For this
end we rewrite these assumptions as a zero finding problem to
which we apply an interval Newton method in {\em the reduced
space}.

Our original phase space is 12 dimensional, the state of bodies is
given by $(q_1, q_2, q_3,\dot{q_1}, \dot{q_2}, \dot{q_3})$. The
center of mass is fixed at the origin. Hence one body's position
and velocity  is determined by other two bodies. We  start from a
collinear position with the third body at the origin and with
equal velocities of the first and the second body (see
\ref{fig:eight2}). Hence it is enough to know  the position and
the velocity of the first body to reconstruct initial condition of
other bodies. Moreover, if we have one solution we could get
another solution by a suitable rotation (both have the same
shape). To remove this degeneracy  we place the first body on the
X axis (we will make computation for rotated Eight, see
Figure~\ref{fig:eight2}). In addition we  fix the size of
trajectory. This fixes the period of the solution. But if we know
one periodic solution then by Kepler law we may obtain a solution
of any period just by rescaling. Finally, we set $q_1(0)=(1,0)$,
hence initial conditions are fixed by the velocity of the first
body.
\par
The {\em reduced space} for Eight is two dimensional and is
parameterized by velocity of second body.  We define map from the
reduced space to full phase space $E:\Real^2 \To \Real^{12}$,
which expands velocity of first body, given by $(v,u)$, to the
initial conditions
$(x_1,y_1,x_2,y_2,x_3,y_3,v_1,u_1,v_2,u_2,v_3,u_3)$ for equation
(\ref{eq:nbp})
\[
E(v,u) = (1,0,-1,0,0,0,v,u,v,u,-2v,-2u).
\]

For each such initial condition exists a solution of the tree body
problem defined on some interval. To each initial configuration,
following that solution, we associate, if it exists, a
configuration in which for the first body the position vector is
orthogonal to its velocity vector for the first time. This defines
the Poincar\'e map $P:\Real^{12} \supset \Omega \To \Real^{12}$.

Now, we define map $R:\Real^{12} \To \Real^2$, by
\[
R(q_1, q_2, q_3, \dot{q}_1, \dot{q}_2, \dot{q}_3) =
(\|q_2-q_1\|^2-\|q_3-q_1\|^2 , (\dot{q}_2-\dot{q}_3)\times q_1),
\]
where by $\times$ we denote vector product, and map $\Phi: \Real^2
\To \Real^2$ by
\[\Phi =R\circ P\circ E.\]
\begin{rem}\label{rem:zero}
If $R(q_1,q_2,q_3,\dot {q}_1,\dot {q}_2,\dot {q}_3)=0$, then $
\|q_2-q_1\|=\|q_3-q_1\|$ and if in addition $y_1=0$ but $x_1\neq
0$, then $\dot{y}_2=\dot{y}_3$. The bodies are then located as in
Fig.~\ref{fig:eight2_e} and conditions (e) and (f) in
Rem.~\ref{rem:eight} are satisfied in a suitably rotated
coordinate frame.
\end{rem}

The following lemma, which is crucial for the proof of
Theorem~\ref{thm:eight} is obtained with computer assistance.
\begin{lem}
\label{lem:zero} There exists  a locally unique $(v,u) \in
\Real^2$ that $\Phi(v,u) = (0,0)$.
\end{lem}
\noindent \textbf{Proof:} We use an interval Newton method
(Theorem~\ref{thm:newton}). First we come close to zero, starting
with some rough initial condition (for example from \cite{S1})
using a  non-rigorous Newton method. Once we have a good candidate
$x_0 = (v_0, u_0)$, we set $[X] = [v_0 - \delta, v_0 + \delta]
\times [u_0 - \delta, u_0 + \delta]$ and  compute rigorously
$\Phi(x_0)$ and $\frac{\partial \Phi([X])}{\partial x}$. For this
purpose we use $C^1$-Lohner algorithm described in \cite{ZLo}. In
this   computation we used a time step $h=0.01$ and the order
$r=7$.

\begin{table}
\begin{tabular}{|c|c|} \hline
 \multicolumn{2}{|c|}{   Result } \\\hline
&\\
  $\Phi(x_0)$ & ([-2.107029e-06, -2.106467e-06],[2.974991e-06, 2.976034e-06]) \\
&\\
  diam $\Phi(x_0)$ & (5.625889e-10, 1.042962e-09) \\
&\\
  $\frac{\partial \Phi([X])}{\partial x}$ &
$\left[\begin{matrix}
[17.622624, 17.643043] & [1.809772,1.827325] \\
[-24.868548,-24.848432] &  [-10.056629,-10.039221]
\end{matrix}\right] $\\

&\\
  $N(x_0,[X])$ & $\begin{matrix}
       ([0.347116886243943,0.347116889993313],\\
       [0.532724941587373,0.532724949187495])\\
  \end{matrix}$ \\
&\\
  diam $N(x_0,[X])$ & (3.749369e-09,7.600121e-09) \\
\hline
\end{tabular}
\vspace*{0.2in} \caption{The data from the proof of
Lemma~\ref{lem:zero}.  \label{tab:8-data}}
\end{table}

It turns out that  the assumption of assertion 1 in
Theorem~\ref{thm:newton} holds for
\newline $x_0 =(0.347116768716, 0.532724944657)$ and $\delta =
10^{-6}$. Numerical data from this computation are listed in
Table~\ref{tab:8-data}.

 Moreover, from Theorem~\ref{thm:newton} we
know that this zero is unique in the set X.\qed

\noindent \textbf{Proof of Theorem \ref{thm:eight}:}  From Lemma
\ref{lem:zero} there exists an initial condition $(v,u)$ in
reduced space that $\Phi(v, u) = (0,0)$. Hence there exists
solution $\bar{q}(t)$ of the three body problem in some interval
$[0,\tilde{T}]$ that in $t=0$ bodies are in collinear position
with third body in the origin and in $t=\tilde{T}$ they are in
isosceles configuration ($r_{12} = r_{13}$). By rotation we could
get solution $\tilde{q}(t)$ that in $t=\tilde{T}$ first body is on
X axis. Obviously we still start in collinear position and end in
isosceles. Now from Remark~\ref{rem:zero} $\tilde{q}(t)$ satisfy
all condition in Lemma~\ref{lem:eight} and hence there exist
$q(t)$ satisfying condition (1),(2),(3).

Local uniqueness follows from local uniqueness in
Lemma~\ref{lem:zero}.\qed

\subsection{Convexity of the Eight}
\begin{thm} \label{thm:convex}
Each lobe of the Eight is convex.
\end{thm}
\noindent \textbf{Proof:} For the proof it is enough to show that
the only inflection point on the curve $q(t)$ is the origin.

 From Lemma~\ref{lem:zero} we get set
$[X]$ which includes initial condition for the Eight in reduced
space. We expand $[X]$ to the full space, we set $\bar{X} =
E([X])$. In the coordinate frame in which the Eight looks like in
Fig.~\ref{fig:eight} and symmetries  $\sigma$ and $\tau$ are
reflections with respect to coordinate axes (see \ref{eq:8-sym})
we see immediately that  the symmetry properties of the Eight
imply that at the origin  $\frac{\partial^2 y_i}{\partial x_i^2} =
0 $ and $\frac{\partial^2 x_i}{\partial y_i^2} = 0 $. To prove
convexity of the Eight we  follow rigorously the trajectory of set
$\bar{X}$ and show that the only point in which $\frac{\partial^2
y_i}{\partial x_i^2} = 0 $ and $\frac{\partial^2 x_i}{\partial
y_i^2} = 0 $ is the origin. The same is true if we rotate the
Eight, as in Figure~\ref{fig:eight2}, to the coordinate system in
which we performed the actual computations and in which we will
work for the remainder of this proof.

Let $q_i(t)=(x_i(t), y_i(t))$ be position of i-th body . If $0
\notin \frac{\partial{x_i}}{\partial t}[t_{k-1},t_{k}]$ then, from
the implicit function theorem, we can write $y_i$ as a function of
$x_i$ for interval between $x_i(t_{k-1})$ and $x_i(t_k)$.
Otherwise we try to represent $x_i$ as a function of $y_i$. For
small enough time steps at least one of this representations is
always possible for the Eight (this is really verified during the
rigorous computations) . Now we derive the formulas for the
derivatives of $y_i(t(x_i))$ with respect to $x_i$ (for
derivatives of $x_i(t(y_i))$ we just exchange $x$ and $y$
variables in formulas).
\begin{eqnarray*}
\frac{\partial y_i(t(x_i))}{\partial x_i} &=& \frac{\partial
y_i}{\partial t} \left(\frac{\partial x_i}{\partial t}\right)^{-1}
\\
 \frac{\partial^2 y_i(t(x_i))}{\partial x_i^2} &=&
\left(\frac{\partial^2 y_i}{\partial t^2} - \frac{\partial^2
x_i}{\partial t^2} \frac{\partial y_i}{\partial x_i} \right)
\left(\frac{\partial x_i}{\partial t}\right)^{-2} \\
\frac{\partial^3 y_i(t(x_i))}{\partial x_i^3}
&=&\frac{\frac{\partial^3 y_i}{\partial t^3}\frac{\partial
x_i}{\partial t} - \frac{\partial^3 x_i}{\partial t^3}
\frac{\partial y_i}{\partial t}+ 2 \left( \frac{\partial^2
x_i}{\partial t^2}\right)^2 \frac{\partial y_i}{\partial x_i} -2
\frac{\partial^2 x_i}{\partial t^2}\frac{\partial^2 y_i}{\partial
t^2} }{ \left(\frac{\partial x_i}{\partial t}\right)^{4}} - \frac{
\left(\frac{\partial^2 x_i}{\partial t^2} \frac{\partial^2
y_i}{\partial x_i^2} \right) }{\left(\frac{\partial x_i}{\partial
t}\right)^{2}}
\end{eqnarray*}
From this equations we need to  know the derivatives $x_i$ and
$y_i$ with respect to time. We obtain them easily by
differentiation of (\ref{eq:nbp}) with respect to time. In fact
this is done during each step  of $C^1$-Lohner algorithm.

By $\varphi(t,x_0)$ we denote the state of bodies (the position
and the velocity) at time $t$ with the initial condition $x_0$ at
$t=0$. Let $h_k$ be length of $k$-th time step, $t_k=h_1 + \dots +
h_k $ - the time after $k$ steps, $[q^k]=
[q_1^k]\times[q_2^k]\times[q_3^k] \times[\dot{q}_1^k]
\times[\dot{q}_2^k]\times[\dot{q}_3^k] \subset \Real^{12}$ be
product of intervals, such that $\varphi(t_k, \bar{X}) \subset
[q^k]$ and $[Q^k]= [Q_1^k]\times[Q_2^k]\times[Q_3^k]
\times[\dot{Q}_1^k] \times[\dot{Q}_2^k]\times[\dot{Q}_3^k] \subset
\Real^{12}$, where $[Q_i^k] = [X_i^k] \times [Y_i^k]$,  be product
of intervals, such that $\varphi([t_{k-1},t_k], \bar{X}) \subset
[Q^k]$ Note that $[Q^k]$ can be seen as an interval enclosure for
whole trajectory between $[q^{k-1}]$ and $[q^k]$, which is
computed during each step of $C^1$-Lohner algorithm.

 For $k$-th step and $i$-th body (except first step for third
body) we check  if at least one of the following conditions is
true
\begin{eqnarray} \label{eq:conv1}
0 \notin \frac{\partial{x_i}}{\partial t}[X_i^k] \mbox{ and } 0
\notin \frac{\partial^2{x_i}}{\partial y_i^2}[X_i^k] \\
 \label{eq:conv2} 0 \notin \frac{\partial{y_i}}{\partial
t}[Y_i^k] \mbox{ and } 0 \notin \frac{\partial^2{y_i}}{\partial
x_i^2}[Y_i^k].
\end{eqnarray}

For first step for third body (starting in the origin) we check if
one of following conditions is satisfied
\begin{eqnarray} \label{eq:conv3}
0 \notin \frac{\partial{x_3}}{\partial t}[X_3^1] \mbox{ and } 0
\in \frac{\partial^2{x_3}}{\partial y_3^2}[X_3^1] \mbox{ and } 0
\notin \frac{\partial^3{x_3}}{\partial y_3^3}[X_3^1] \\
 \label{eq:conv4} 0 \notin \frac{\partial{y_3}}{\partial
t}[Y_3^1] \mbox{ and } 0 \in \frac{\partial^2{y_3}}{\partial
x_3^2}[Y_3^1] \mbox{ and } 0 \notin
\frac{\partial^3{y_3}}{\partial x_3^3}[Y_3^1].
\end{eqnarray}
To verify above conditions we  follow rigorously using
$C^1$-Lohner algorithm  the trajectory of set $\bar{X}$ until
reaching the section described in proof of Lemma~\ref{lem:zero}
and for each time step and each body we check suitable condition.
It turns out that for each time step and each body these
conditions were satisfied. This finishes the proof. \qed

\subsection{Some numerical data from the convexity of Eight computation}
The parameters of the methods were: a time step $h=0.01$, order
$r=7$. We needed 53 steps to cross the section, below we show data
for  characteristic cases in order to show that it was really
quite easy to verify, with  algorithms we used.

 \noindent\textbf{Step 1.} We start in
collinear position with third body in the origin. So we are in an
{\em inflection point}.
\begin{table}
\begin{center}
\begin{tabular}{|c|c|c|c|} \hline
  $i$ & $\frac{\partial{x_i}}{\partial t}[X_i^1] $&
   $\frac{\partial^2{x_i}}{\partial y_i^2}[X_i^1]$ &   $ \frac{\partial^3{x_i}}{\partial y_i^3}[X_i^1]$
   \\\hline
  1 & [0.334402,0.347118] & [15.3592,17.9897] & [136.616,219.114] \\
  2 & [0.347116,0.360049] & [-16.5013,-14.111] & [119.951,192.562] \\
  3 & [-0.695034,-0.69385] & [-0.0682713,0.269952] & [-30.969,-26.3718]\\ \hline
\end{tabular}
\vspace*{0.2in} \caption{The data from the proof of
Theorem~\ref{thm:convex}, for step 1, neighborhood of deflection
point. \label{tab:convex-step1}}
\end{center}
\end{table}
Numerical data for this case are given in
Table~\ref{tab:convex-step1}. We see that for third body we have
$0 \in \frac{\partial^2{x_i}}{\partial y_i^2}[X_i^1]$, but this
derivative is monotonic ($ \frac{\partial^3{x_i}}{\partial
y_i^3}[X_i^1] <0 $), hence there can be only one zero of it in
interval $[X_3^1]$. But we know that one zero is at $(0,0)\in
[X_3^1]$.

\noindent \textbf{Steps 2-36 and 38-53} In this case all second
derivatives do not contain 0. Table~\ref{tab:convex-step2}
contains data obtained in second step.
\begin{table}
\begin{center}
\begin{tabular}{|c|c|c|c|} \hline
  $i$ & $\frac{\partial{x_i}}{\partial t}[X_i^2]$ &
   $\frac{\partial^2{x_i}}{\partial y_i^2}[X_i^2]$ &    $\frac{\partial^3{x_i}}{\partial y_i^3}[X_i^2]$
   \\ \hline
  1 & [0.32222,0.334722]& [16.7085,19.6225] & [155.007,250.021] \\
  2 & [0.35972,0.372882] & [-15.2203,-13.0177] & [105.778,171.191] \\
  3 & [-0.695669,-0.69444] & [0.126359,0.472046] & [-30.9533,-26.1203] \\ \hline
\end{tabular}
\vspace*{0.2in} \caption{The data from the proof of
Theorem~\ref{thm:convex}, for step 2.\label{tab:convex-step2}}
\end{center}
\end{table}

\noindent \textbf{Step 37.} In this case we cannot represent $y_1$
as function of $x_1$, hence we interchange variables and represent
$x_1$ as function of $y_1$ and  check condition \ref{eq:conv2}
instead of \ref{eq:conv1}. The first body is in the rightmost
position. Table~\ref{tab:convex-step37} contains the derivatives
for this case.
\begin{table}
\begin{center}
\begin{tabular}{|c|c|c|c|} \hline
  $i$ & $\frac{\partial{x_i}}{\partial t}[X_i^2]$ &
   $\frac{\partial^2{x_i}}{\partial y_i^2}[X_i^2]$ &    $\frac{\partial^3{x_i}}{\partial y_i^3}[X_i^2]$
   \\ \hline
  1 & [-0.00287209,0.00468403]& - & - \\
  2 & [0.904939,0.922079] & [-2.55715,-2.16831] & [4.35197,10.6201] \\
  3 & [-0.919428,-0.909617] & [2.56371,3.03259] & [-3.51203,3.48564] \\ \hline
  $i$ & $\frac{\partial{y_i}}{\partial t}[Y_i^2]$ &
  $\frac{\partial^2{y_i}}{\partial x_i^2}[Y_i^2]$ &    $\frac{\partial^3{y_i}}{\partial
  x_i^3}[Y_i^2]$\\ \hline
  1 & [0.480975,0.48288] & [-3.24824,-3.11737] & [-2.98453,-0.860616] \\ \hline
\end{tabular}
\vspace*{0.2in} \caption{The data from the proof of
Theorem~\ref{thm:convex}, for step 37.\label{tab:convex-step37}}
\end{center}
\end{table}

\section{Doubly symmetric choreographies with even number of bodies}
\label{sec:linchains}
\subsection{Symmetries} Many of the choreographies found by Simo
\cite{S1} have at least one symmetry. Just as in case of the Eight
this is not only a symmetry of trajectory image, but also of how
the curve describing the trajectory is parameterized with time.

In this section we introduce a notation for symmetries which will
be used till the end of this paper. By $S_x$, $S_y$, $S_0$ we want
to denote a symmetry with respect to the x axis, the y axis and to
the origin. These spatial symmetries act also on time variable
parameterizing curves. \newline
 Let $T$ be any positive
real number. We define actions of $S_x$, $S_y$ and $S_0$ on
$\Real/T\Integer$ and on $\Real^2$ as follows: \DM
\begin{array}{ll}
S_x(t) = -t+\frac{T}{2}, & S_x(x, y) = (-x, y),\\
S_y(t) = -t, & S_y(x, y) = (x,-y),\\
S_0(t) = t+\frac{T}{2}, & S_0(x, y)=(-x, -y).
\end{array}
\EDM It follows from this definition that $S_0 = S_x \circ S_y =
S_y \circ S_x$.

Let $q(t): \Real/T\Integer \To \Real^2$. We say that $q(t)$ is
invariant (equivariant) with respect to the action of S if
$S(q(t))=q(S(t))$ for all $t$. Moreover we assume that $q(t)$ is
$C^2$ function. If $q(t)$ is invariant with respect to $S_x$
(resp. $S_y$) then $\dot{q}(S_x(t)) = S_y(\dot{q}(t))$ (resp.
$\dot{q}(S_y(t)) = S_x(\dot{q}(t))$). Hence $S_0$ invariance
implies that $\dot{q}(S_0(t)) = S_0(\dot{q}(t))$.

From now on we will enumerate bodies starting from 0. We set
$q_i(t)=q\left(t+\frac{T}{N} \cdot i \right)$ for $i = 0, \dots,
N-1$.

\subsection{Doubly symmetric choreographies with even number of bodies}
\label{subsec:linchains}

We will consider only cases with even number of bodies.

 Let $T= N \cdot \bar{T}>0$. We search for a function $q(t): \Real/T\Integer \To \Real^2$ which
has following properties:
\par\noindent
\begin{enumerate}
\item[P1.]
 for each $t$ the origin is center of mass,
\BE \sum_{i=0}^{N-1}{q\left(t+ \bar{T} \cdot i \right)}=0, \EE

\item[P2.] $q(t)$ is invariant with respect to
\begin{enumerate}
\item $S_x$ i.e. $q(S_x(t)) = S_x(q(t))$,
\item $S_0$ i.e. $q(S_0(t)) = S_0(q(t))$,
\end{enumerate}
\item[P3.] function $x=(q_0(t), q_2(t),\dots, q_{N-1})$ where
$q_i(t)=q\left(t+\bar{T} \cdot i \right)$ for $i = 0, \dots, N-1$,
is a T-periodic solution of the Newtonian N-body problem
(\ref{eq:nbp}).
\end{enumerate}

Observe that condition \textbf{P3} says that bodies trace each
other with a constant phase shift.

\begin{figure}
 \centerline{\includegraphics[scale = 0.7]{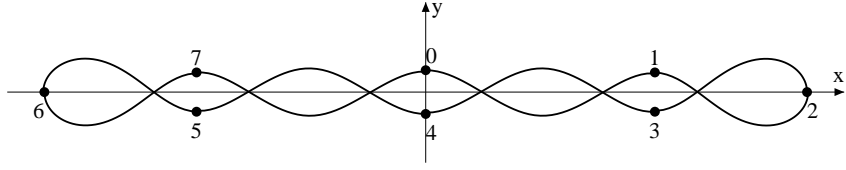}}
 \caption{Linear chain with 8 bodies. \label{fig:8body}}
\end{figure}

The following Lemma, which is analogous to Lemma~\ref{lem:eight},
gives necessary and sufficient conditions for an existence of a
choreography satisfying P1,P2 and P3.
\begin{lem}
\label{lem:linchains} Let $N=2n$ be a number of bodies. There
exists a function q(t) with properties P1, P2, P3 if and only if
there are functions $q_i:[0, \bar{T}/2]\To R^2$ for
$i=0,1,\dots,N-1$ such as:
\begin{enumerate}
\item
$q_0(0)=(0,y_0)$ for some $y_0 \neq 0$ ($q_0(0)$ is on Y axis),
\item
At time $t_0=0$ for $i=0,1,\dots,N/2-1$ we have
\begin{enumerate}
\item $q_i(t_0)=S_x(q_{\frac{N}{2}-i}(t_0))$
\item $\dot{q}_i(t_0)=S_y(\dot{q}_{\frac{N}{2}-i}(t_0))$
\item $q_i(t_0) = S_0(q_{\frac{N}{2}+i}(t_0))$
\item $\dot{q}_i(t_0) = S_0(\dot{q}_{\frac{N}{2}+i}(t_0))$
\end{enumerate}

\item
At time $t_1=\bar{T}/2$ for $i=0,1,\dots,N/2-1$ we have
\begin{enumerate}
\item $q_i(t_1)=S_x(q_{\frac{N}{2}-i-1}(t_1))$
\item $\dot{q}_i(t_1)=S_y(\dot{q}_{\frac{N}{2}-i-1}(t_1))$
\item $q_i(t_1) = S_0(q_{\frac{N}{2}+i}(t_1))$
\item $\dot{q}_i(t_1) = S_0(\dot{q}_{\frac{N}{2}+i}(t_1))$
\end{enumerate}
\item $x(t) = (q_0(t),q_1(t),\dots,q_{N-1}(t))$ is a solution of the Newton N-body problem for $t \in [0, \bar{T}/2]$.
\end{enumerate}
\end{lem}

\noindent \textbf{Proof:} If we have a trajectory with properties
P1, P2, P3 then it is easy to show that functions $q_i(t)$ defined
in P3 satisfy conditions (1)-(3). On the other hand if we have a
functions $q_i(t)$ which have properties (1)-(3) then we define
$q(t)$ by
\begin{equation}
q(t)=\left\{ \begin{array}{ll} q_i(t-i\bar{T}) & \mbox{for } t\in
[
i\bar{T}, (i+\frac{1}{2})\bar{T}] \mbox{ and } i = 0, 1, \dots, N-1\\
S_x(q_{\frac{N}{2}-i}(i\bar{T}-t))& \mbox{for } t\in [
(i-\frac{1}{2})\bar{T}, i\bar{T}] \mbox{ and } i = 1, 2, \dots, \frac{N}{2}\\
S_x(q_{\frac{3N}{2}-i}(i\bar{T}-t))& \mbox{for } t\in [
(i-\frac{1}{2})\bar{T}, i\bar{T}] \mbox{ and } i = \frac{N}{2}+1,
\dots, N
\end{array}\right.
\end{equation}\qed

From Lemma~\ref{lem:linchains} it follows that the proof of an
existence of a doubly symmetric choreography is equivalent  to
some boundary value problem for the N-body equation. We will now
formulate this problem as a zero finding problem for a suitable
map.

Original phase space for the planar N-body problem has $4N$
dimensions. It turns out that, if initial conditions satisfy all
conditions in point 2 of Lemma~\ref{lem:linchains} it is enough to
know values of only $N$ variables to recover rest of them. We
still may obtain solutions of any period, to determine this we can
fix the size of curve by fixing one variable. Hence our {\em
reduced space} be $(N-1)$-dimensional. In the next paragraph we
will be more specific.

We define map $E:\Real^{N-1} \To \Real^{4N}$, which expands, using
symmetries from (2), initial conditions from reduced space to the
full phase space: \newline $(x_0,y_0,\dot{x}_0,\dot{y}_0,
x_1,y_1,\dot{x}_1,\dot{y}_1, \dots,
x_{N-1},y_{N-1},\dot{x}_{N-1},\dot{y}_{N-1})$. We consider two
cases: $N=4k$ and $N=4k+2$.

\par\noindent For $N=4k$ we set
\begin{eqnarray*}
 E(\dot{x}_0 \times \prod_{i=1}^{k-1}
(x_i,y_i,\dot{x}_i,\dot{y}_i) \times (x_k,\dot{y}_k))
 =
(0, a, \dot{x}_0, 0) \times \prod_{i=1}^{k-1}
(x_i,y_i,\dot{x}_i,\dot{y}_i) \\
\times (x_k, 0, 0, \dot{y}_k) \times \prod_{i=1}^{k-1}
(x_{k-i},-y_{k-i},-\dot{x}_{k-i},\dot{y}_{k-i})\\
\times (0, -a, -\dot{x}_0, 0)  \times \prod_{i=1}^{k-1}
(-x_i,-y_i,-\dot{x}_i,-\dot{y}_i)\\
\times (-x_k, 0, 0, - \dot{y}_k) \times \prod_{i=1}^{k-1}
(-x_{k-i},y_{k-i},\dot{x}_{k-i},-\dot{y}_{k-i})
\end{eqnarray*}

\par\noindent For $N=4k+2$ we set
\begin{eqnarray*}
E(\dot{x}_0 \times \prod_{i=1}^k (x_i,y_i,\dot{x}_i,\dot{y}_i))
=(0, a, \dot{x}_0, 0) \times \prod_{i=1}^k
(x_i,y_i,\dot{x}_i,\dot{y}_i)\\
\times \prod_{i=0}^{k-1}
(x_{k-i},-y_{k-i},-\dot{x}_{k-i},\dot{y}_{k-i}) \times (0, -a,
-\dot{x}_0, 0)  \\\times \prod_{i=1}^k
(-x_i,-y_i,-\dot{x}_i,-\dot{y}_i)  \times \prod_{i=0}^{k-1}
(-x_{k-i},y_{k-i},\dot{x}_{k-i},-\dot{y}_{k-i})
\end{eqnarray*}
In both cases $a$ is a parameter fixing a size of an orbit.

Again, like for the Eight, for every such initial condition exists
a solution of the $N$-body problem (\ref{eq:nbp}) defined on some
interval. To every initial configuration, following that solution,
we associate, if it exists, a configuration in which for the first
time
\begin{itemize}
\item for $N=4k$: bodies $k$ and $k-1$ have equal $x$ coordinate
($x_k=x_{k-1}$),
\item for $N=4k+2$: $k$-th body is on the X axis ($y_k=0$).
\end{itemize}
This defines the Poincar\'e map $P:\Real^{4N} \supset \Omega \To
\Real^{4N}$.

Now, we define the reduction map $R:\Real^{4N} \To \Real^{N-1}$ in
such way that $R$ has zeroes in points satisfying conditions (3a)
and (3b) in Lemma~\ref{lem:linchains} and only in such points.
Observe that we don't have to worry about conditions  (3c) and
(3d) in Lemma~\ref{lem:linchains}, because from the properties of
(\ref{eq:nbp}) if follows that if (2c) and (2d) holds, then (3c)
and (3d) are satisfied for {\em any} $t_1$.

\noindent For $N= 4k$ we set
\begin{eqnarray*}
R\left(\prod_{i=0}^{N-1}(x_i,y_i,\dot{x}_i,\dot{y}_i)\right)
 =( y_k + y_{k-1}, \dot{x}_k + \dot{x}_{k-1},\dot{y}_k - \dot{y}_{k-1})\\
 \times \prod_{i=0}^{k-2}(x_i - x_{2k-i-1}, y_i + y_{2k-i-1}, \dot{x}_i + \dot{x}_{2k-i-1},\dot{y}_i - \dot{y}_{2k-i-1})
\end{eqnarray*}

\noindent For $N= 4k+2$ we set
\[
R\left(\prod_{i=0}^{N-1}(x_i,y_i,\dot{x}_i,\dot{y}_i)\right)
 = \{v_k\}\times \prod_{i=0}^{k-1}(x_i - x_{2k-i}, y_i + y_{2k-i}, \dot{x}_i + \dot{x}_{2k-i},\dot{y}_i - \dot{y}_{2k-i})
\]
We define map $\Phi: \Real^{N-1} \supset E^{-1}(\Omega) \To
\Real^{N-1}$ by
\[\Phi =R\circ P\circ E.\]

\begin{thm}\label{thm:linzero}
If for some  $x \in \Real^{N-1}$ $\Phi(x)=0$, then there exists a
trajectory with properties P1, P2, P3.
\end{thm}
\noindent \textbf{Proof:} Let $x$ be point in which we have
$\Phi(x)=0$. Then $E(x)$ satisfies condition (1),(2) in
Lemma~\ref{lem:linchains} and from construction of $\Phi$ we
obtain that there exist trajectories ($q_i(t)$) starting with
initial condition $E(x)$ and ending with $\hat{x}$ such that
$R(\hat{x})=0$, hence  condition (3) in Lemma~\ref{lem:linchains}
is satisfied. Now the assertion follows from
Lemma~\ref{lem:linchains}.  \qed

\section{Existence of the SuperEight - Gerver orbit}
\label{sec:gerver} 
\begin{figure}[hdtb]
 \centerline{\includegraphics[scale = 0.6]{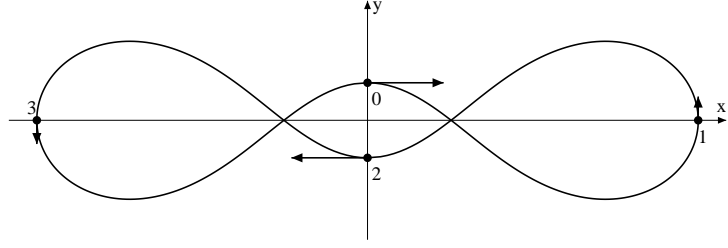}}
 \caption{Gerver orbit - initial position \label{fig:gerver}}
\end{figure}

\begin{figure}
 \centerline{\includegraphics[scale = 0.6]{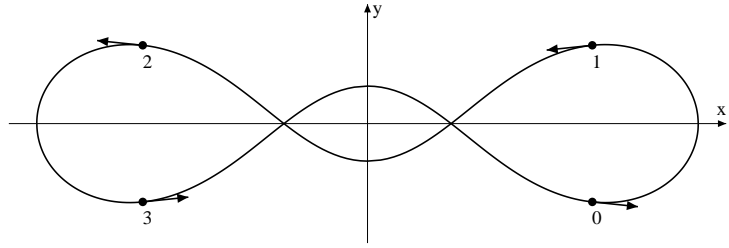}}
 \caption{Gerver orbit - final position \label{fig:gerver_e}}
\end{figure}

The Gerver orbit (Figure \ref{fig:gerver},
\ref{fig:gerver_e}) is a choreography with 4 bodies forming a
linear chain. It's the simplest trajectory after the Eight.

With a computer assistance we proved the following
\begin{thm}
\label{thm:gerver}
  The Gerver SuperEight exists and is locally unique (up to obvious
  symmetries and rescaling).
\end{thm}

We show the existence of SuperEight using an approach described in
Section~\ref{sec:linchains}, i.e. we verify assumptions of
Theorem~\ref{thm:linzero}. Below we give some details.

We set
\begin{equation}\label{eq:Egerver}
  E(x_1,\dot{x}_0,\dot{y}_1)=(0,a,\dot{x}_0,0,
  x_1,0,0,\dot{y}_1,
  0,-a,-\dot{x}_0,0,
  -x_1,0,0,-\dot{y}_1),
\end{equation}
where $a$ is  a parameter fixing the size of the orbit.
\begin{equation}\label{eq:Rgerver}
 R(x_0,y_0,\dot{x}_0,\dot{y}_0,x_1,y_1,\dot{x}_1,\dot{y}_1,x_2,y_2,\dot{x}_2,\dot{y}_2,x_3,y_3,\dot{x}_3,\dot{y}_3)=(y_1+y_0,\dot{x}_1+\dot{x}_0,\dot{y}_1-\dot{y}_0),
\end{equation}
Poinc\'{a}re section $S$ is defined by
\begin{equation}\label{eq:Sgerver}
 S(x_0,y_0,\dot{x}_0,\dot{y}_0,x_1,y_1,\dot{x}_1,\dot{y}_1,x_2,y_2,\dot{x}_2,\dot{y}_2,x_3,y_3,\dot{x}_3,\dot{y}_3)=x_1-x_0=0.
\end{equation}

 At first we have proved an existence of a zero of $\Phi(x)$ using an interval
Newton method, but in the paper we present data from proof based
on Krawczyk method. In case of Gerver orbit the choice of the
method isn't so important, because if we take a time step small
enough or smaller set $[X]$, then the computed matrix
$\left[\frac{\partial \Phi}{\partial x}([X])\right]$ becomes
invertible and the proof usually goes through. But in general the
main problem of application of an interval Newton method is that
of an invertibility of $\left[\frac{\partial \Phi}{\partial
x}([X])\right]$. We avoid this using the Krawczyk method.

To check the assertion (2) in Theorem~\ref{thm:Krawczyk} we used
the $C^1$-Lohner algorithm with  order $r=6$. The time step was
set to $h=0.002$. As matrix $C$ we used an inverse of the
monodromy matrix computed nonrigorously in a point $\bar{x}$, i.e.
$C = {\frac{\partial \Phi (\bar{x})}{\partial x}}^{-1}$.
\begin{table}
\begin{center}
\begin{tabular}{|c|l|} \hline
 \multicolumn{2}{|c|}{Initial values} \\ \hline
 &\\
 $\bar{x}$ &  (1.382857, 1.87193510824,0.584872579881)\\
 $a$ &  0.157029944461 \\
 $[X]$ & $ \bar{x} + [-10^{-7},10^{-7}]^3$ \\
  \hline
\end{tabular}
\end{center}
\vspace*{0.2in}\caption{Data from the proof of existence of Gerver
SuperEight. Initial values.   \label{tab:gerverd1} }
\end{table}

\begin{table}
\begin{center}
\begin{tabular}{|c|l|} \hline
 \multicolumn{2}{|c|}{ Computed values } \\ \hline
&\\
  $C$ & $\left[\begin{matrix}
     -2.15400 & 0.257911 & 0.786925 \\
     -0.08163 & 0.293713 & 0.043565 \\
     0.939059 & -0.10027 & 0.158399
    \end{matrix}\right]$ \\
&\\

  $\Phi(\bar{x})$ & $
   \left[\begin{matrix}
[-2.87020e-09,-2.26613e-09]\\
[-1.21155e-08,-1.06812e-08]\\
[-5.45542e-08,-5.10016e-08]
   \end{matrix} \right]$ \\
&\\
 diam $\Phi(\bar{x})$ &   $\left[ \begin{matrix}
 6.04064e-10\\1.43432e-09\\3.55268e-09
  \end{matrix} \right]$ \\
&\\
 $\frac{\partial \Phi}{\partial x}([X])$ &
$\left[\begin{matrix}
[-0.1664,-0.1657] & [0.39070,0.39119] & [0.71771,0.71790] \\
[-0.1764,-0.1750] & [3.52548,3.52654] & [-0.0968,-0.0964] \\
[0.87189,0.87534] & [-0.0867,-0.0842] & [1.99599,1.99697]
\end{matrix}\right]$ \\
&\\
  $K(\bar{x},[X],\Phi)$  &
  $\left[ \begin{matrix}
  [1.382857036247056692,1.382857041633411832]\\
  [1.871935113301492981,1.871935114053588922]\\
  [0.5848725887384301769,0.5848725902808686872]
  \end{matrix} \right]$ \\
  &\\
 diam $K(\bar{x},[X],\Phi)$ &
 $\left[ \begin{matrix}
 5.386355139691545446e-09\\
 7.520959410811656198e-10\\
 1.54243851024915557e-09
  \end{matrix} \right]$ \\
  \hline
\end{tabular}
\end{center}
\vspace*{0.2in} \caption{Data from the proof of existence of
Gerver SuperEight. Matrix $C$ and results of computation.}
\label{tab:gerverd2}
\end{table}
Tables \ref{tab:gerverd1} and \ref{tab:gerverd2} contain numerical
data from this proof.

\section{Existence of the 'Linear chain' orbit for the 6 bodies}
\label{sec:6body} 
\begin{figure}[hdtb]
 \centerline{\includegraphics[scale = 0.6]{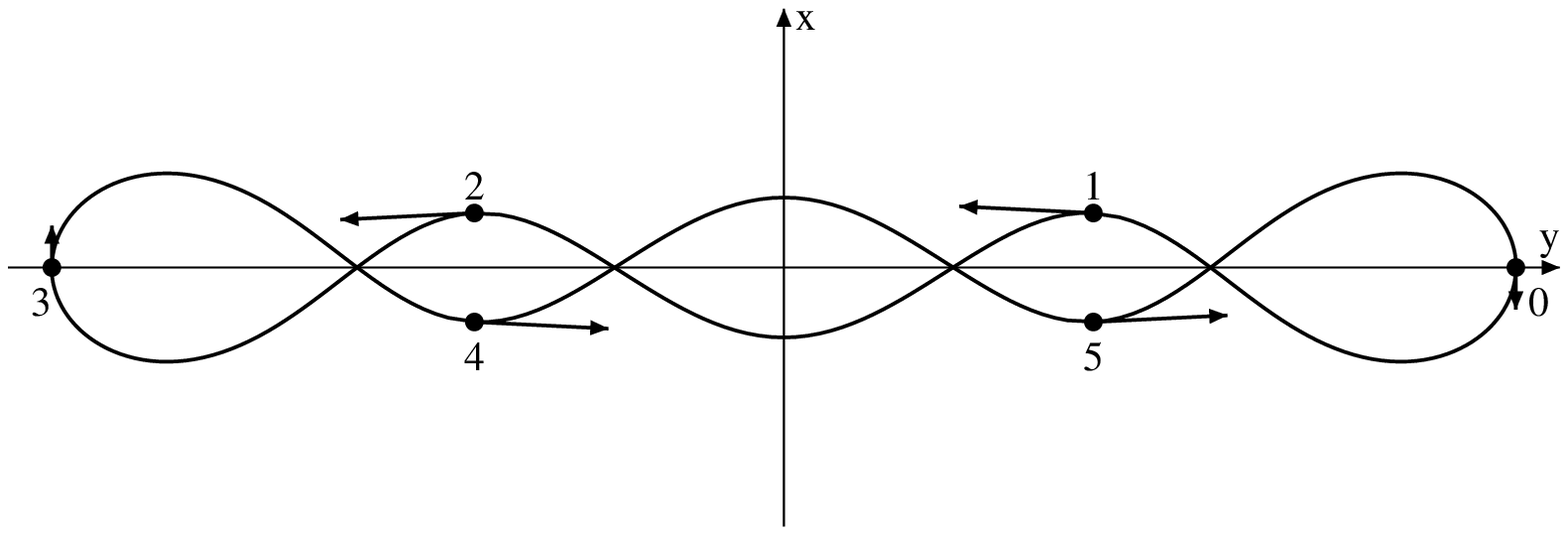}}
 \caption{'Linear chain' orbit for the 6 bodies - initial position \label{fig:6body}}
\end{figure}

Figure \ref{fig:6body}  displays a  linear chain
choreography with 6 bodies.

\begin{figure}
 \centerline{\includegraphics[scale = 0.6]{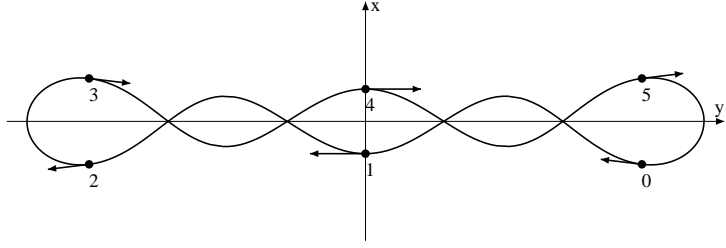}}
 \caption{'Linear chain' orbit for the 6 bodies - final position \label{fig:6body_e}}
\end{figure}

In this section we report about the computer assisted proof of the
following
\begin{thm}
\label{thm:cl6b}
  The linear chain for $6$ bodies exists and is locally unique (up to obvious
  symmetries and rescaling).
\end{thm}

Again, we show an existence of this orbit using an approach
described in Section~\ref{sec:linchains} with some minor changes.
To speed up calculation and to increase accuracy we take into
account in computation that all time $q_3(t)=-q_0(t)$,
$q_4(t)=-q_1(t)$ and $q_5(t)=-q_2(t)$ (we change equation
(\ref{eq:nbp}) doing a suitable substitution). So, our full space
for ODE has now 12 dimensions. We use also a different time
parameterization (we do time shift of $\frac{1}{4}$ of the
period). To use exactly approach described in
Section~\ref{sec:linchains} we interchange axes (see Fig.
\ref{fig:6body} and \ref{fig:6body_e}). From
Lemma~\ref{lem:linchains} we will get, in this coordinate frame, a
doubly symmetric periodic solution $q(t)$. It is easy to see that
$\bar{q} = q(t - \frac{T}{4})$ is then a solution sharing needed
symmetries in the original coordinate frame. All data are given in
frame with interchanged axes.

 We set
\begin{equation}\label{eq:E6body}
  E(\dot{x}_0,x_1,y_1,\dot{x}_1,\dot{y}_1) =
  (0,a,\dot{x}_0,0,x_1,y_1,\dot{x}_1,\dot{y}_1,x_1,-y_1,-\dot{x}_1,\dot{y}_1),
\end{equation}
where $a$ is a parameter fixing the size of the orbit.
\begin{equation}\label{eq:R6body}
 R(x_0,y_0,\dot{x}_0,\dot{y}_0,x_1,y_1,\dot{x}_1,\dot{y}_1,x_2,y_2,\dot{x}_2,\dot{y}_2)
 =(\dot{x}_1,x_0-x_2,y_0+y_2,\dot{x}_0+\dot{x}_2,\dot{y}_0-\dot{y}_2).
\end{equation}
Poinc\'{a}re section $S$ is defined by
\begin{equation}\label{eq:S6body}
 S(x_0,y_0,\dot{x}_0,\dot{y}_0,x_1,y_1,\dot{x}_1,\dot{y}_1,x_2,y_2,\dot{x}_2,\dot{y}_2)=y_1=0.
\end{equation}

 To find a zero of $\Phi(x)$ we use Krawczyk method. We check the assertion (2)
 in Theorem~\ref{thm:Krawczyk} using $C^1$-Lohner algorithm
 \cite{ZLo} of order $r=9$ and a time step $h=0.0025$ for a computation in point $\bar{x}$ and
 $h=0.001$ for a computation on set $[X]$. As matrix $C$ we used an inverse of the
monodromy matrix computed nonrigorously in a point $\bar{x}$, i.e.
$C = {\frac{\partial \Phi (\bar{x})}{\partial x}}^{-1}$.

In Tables \ref{tab:c6bd1} and \ref{tab:c6bd2} we present data from
computation of Krawczyk method.
\begin{table}
\begin{center}
\begin{tabular}{|c|l|}\hline
 & Initial value \\ \hline
 &\\
 $\bar{x}$ &  $\left[\begin{matrix}
 -0.635277524319\\
 0.140342838651\\
 0.797833002006\\
 0.100637737317\\
 -2.03152227864

 \end{matrix} \right]$ \\
 &\\
 $a$ &  1.887041548253914 \\
 &\\
 $[X]$ &  $\left[\begin{matrix}
[-0.635277525319,-0.635277523319]\\
[0.140342837651,0.140342839651]\\
[0.797833001006,0.797833003006]\\
[0.100637736317,0.100637738317]\\
[-2.03152227964,-2.03152227764]

 \end{matrix} \right]$ \\
 &\\
 diam $[X]$ & $\left[\begin{matrix}
 2.0e-09\\
 2.0e-09\\
 2.0e-09\\
 2.0e-09\\
 2.0e-09
 \end{matrix} \right]$ \\
\hline
\end{tabular}
\end{center}
\vspace*{0.2in} \caption{Data from the proof of existence of
linear chain for 6 bodies. Initial values for Krawczyk method.}
\label{tab:c6bd1}
\end{table}

\begin{table}
\begin{tabular}{|c|l|}\hline
 & Computed value \\ \hline
 &\\
  $\Phi(\bar{x})$ & $
   \left[\begin{matrix}
  [-3.1311957909658e-11,3.156277062700585e-11]\\
  [-4.528821762050939e-12,4.574757239694804e-12]\\
  [-1.063704679893362e-11,1.051470022161993e-11]\\
  [-3.084105193451592e-11,3.117495150917193e-11]\\
  [-1.203726007759087e-11,1.193112275643671e-11]
   \end{matrix} \right]$ \\ 
   &\\
 $\diam \Phi(\bar{x})$ &   $\left[ \begin{matrix}
  6.287472853666386e-11\\
  9.103579001745743e-12\\
  2.115174702055356e-11\\
  6.201600344368785e-11\\
  2.396838283402758e-11
  \end{matrix} \right]$
  \\ 
  &\\
  $K(\bar{x},[X],\Phi)$  &
  $\left[ \begin{matrix}
  [-0.6352775243616679557,-0.6352775242763283314]\\
  [0.1403428386430521646,0.1403428386590999943]\\
  [0.797833001999263769,0.797833002012834469]\\
  [0.10063773728817425324,0.1006377373457752189]\\
  [-2.031522278710178764,-2.031522278575771612]
  \end{matrix} \right]$ \\ 
   & \\
 diam $K(\bar{x},[X],\Phi)$ &
 $\left[ \begin{matrix}
 8.53396242561643703e-11\\
 1.604782973174678773e-11\\
 1.357070011920313846e-11\\
 5.760096566387318262e-11\\
 1.344071520748002513e-10
  \end{matrix} \right]$ \\
\hline
\end{tabular}
\vspace*{0.2in} \caption{Data from the proof of an existence of
linear chain for 6 bodies. Results of computation of Krawczyk
method.} \label{tab:c6bd2}
\end{table}

\section{Some technical data}
\label{sec:tech-data} All computations were performed on {\em AMD
Athlon 1700XP } with 256 MB DDRAM memory, with Windows 98SE
operating system. We used CAPD package\cite{Capd} and  Borland C++
5.02 compiler.

In the listing below $r$ is an order and $h$ is a time step used
in the $C^1$-Lohner algorithm \cite{ZLo}.

The computation times for The Eight, $h=0.01$, $r=7$
\begin{itemize}
\item in point $\bar{x}$ : 1.417 sec
\item  for set $[X]$ : 2.66 sec
\item convexity : 1.15501 sec
\end{itemize}

For the proof of an existence of Gerver solution in 4-body problem
we used  $r=6$, $h=0.002$. The computation times for both
$\bar{x}$ and set $[X]$ were approximately equal to $30.5$
seconds.

For the proof of linear chain for 6-body problem
\begin{itemize}
\item computation
of Poincar\'e Map for set $[X]$ took $57.5$ seconds with $h=0.001$
and $r=9$
\item computation for $\bar{x}$ took $23.8$
seconds with $h=0.0025$ and $r=9$
\end{itemize}

The programm performing the proofs is available on {\em
http://www.ap.krakow.pl/\~{}tkapela}



\begin{thebibliography}{MMMM}

\bibitem[A]{A} G. Alefeld, \emph{ Inclusion methods for systems of nonlinear equations - the interval
   Newton method and modifications.\/} in \emph{ Topics in Validated Computations}
   J. Herzberger (Editor), 1994 Elsevier Science B.V.


\bibitem[Capd]{Capd} CAPD - Computer assisted proofs in dynamics, a package for rigorous numerics,
        {\em http://limba.ii.uj.edu.pl/\~{}capd}


\bibitem[Ch]{Ch} A. Chenciner, \emph{Some facts and more questions about the "Eight"}, Proceedings
of the conférence "Nonlinear functional analysis", Taiyuan  2002,
(World Scientific, in press)


\bibitem[CGMS]{CGMS} A. Chenciner, J. Gerver, R. Montgomery, C. Sim\'o \emph{Simple Choreographic Motions of
$N$ Bodies: A Preliminary Study},. Geometry, mechanics, and
dynamics, 287--308, Springer, New York, 2002.


\bibitem[CM]{CM} A. Chenciner and R. Montgomery, \emph{A remarkable periodic solution of
 three-body problem in the case of equal masses}, Annals of
 Mathematics, 152 (2000),881--901

\bibitem[G1]{Gp1}Z. Galias. \emph{Interval methods for rigorous investigations of
periodic orbits,}  Int. J. Bifurcation and Chaos, 11(9):2427-2450,
2001

\bibitem[G2]{Gp2} Z. Galias, \emph{Rigorous investigations of Ikeda map by means of interval
arithmetic,}  Nonlinearity, 15:1759-1779, 2002


\bibitem[K]{K} R. Krawczyk, \emph{Newton-Algorithmen zur Bestimmung von
Nullstellen mit Fehlerschanken}, Computing 4, 187--201 (1969)

\bibitem[M]{M} C. Moore, \emph{Braids in Classical Gravity},
Physical Review Letters, \textbf{70} (1993), 3675--3679

\bibitem[Mo]{Mo} R.E. Moore, {\em Interval Analysis.} Prentice
Hall, Englewood Cliffs, N.J., 1966

\bibitem[MR]{MR} R. Montgomery,  \emph{ A new solution to the three-body
problem.} Notices Amer. Math. Soc. 48 (2001), no. 5, 471--481.

\bibitem[N]{N} A. Neumeier, Interval methods for systems of
equations. Cambrigde University Press, 1990.

\bibitem[S]{S} C. Sim\'o, \emph{Choreographies of the N-body problem,}
 http://www.maia.ub.es/dsg/nbody.html

\bibitem[S1]{S1}
C. Sim\'o, \emph{Periodic orbits of the planar N-body problem with
equal masses and all bodies on the same path},  in {\em The
Restless Universe: Applications of N-Body Gravitational Dynamics
to Planetary, Stellar and Galactic Systems},265--284, ed. B.
Steves and A. Maciejewski, NATO Advanced Study Institute, IOP
Publishing, Bristol, 2001, see also
http://www.maia.ub.es/dsg/2001/


\bibitem[S2]{S2}C. Sim\'o,  \emph{New families of solutions in $N$-body problems.} European
Congress of Mathematics, Vol. I (Barcelona, 2000), 101--115,
Progr. Math., 201, Birkhäuser, Basel, 2001.

\bibitem[S3]{S3} C. Sim\'o, \emph{Dynamical properties of the figure eight solution of
the three-body problem.} Celestial mechanics (Evanston, IL, 1999),
209--228, Contemp. Math., 292, Amer. Math. Soc., Providence, RI,
2002.


\bibitem[ZLo]{ZLo}  P.\ Zgliczy\'nski,
\emph{$C^1$-Lohner algorithm},  Foundations of Computational
Mathematics,  (2002) 2:429--465

\end{thebibliography}
\end{document}